\documentclass[11pt]{article}
\usepackage{enumerate}
\usepackage{amssymb,a4wide,latexsym,makeidx,epsfig,fleqn}
\usepackage{amsthm}
\usepackage{amsmath}
\usepackage{enumerate}
\newtheorem{theorem}{Theorem}[section]
\newtheorem{remark}[theorem]{Remark}

\newtheorem{definition}[theorem]{Definition}
\newtheorem{lemma}[theorem]{Lemma}

\begin{document}
\textwidth 150mm \textheight 225mm
\title{Inertia indices of a complex unit gain graph in terms of matching number  \footnote{This work is supported by the National Natural Science Foundations of China (No. 11901253), the Natural Science Foundation for Colleges and Universities in Jiangsu Province of China (No. 19KJB110009), and the Science Foundation of Jiangsu Normal University (No.18XLRX021).}}
\author{{Yong Lu and Qi Wu\footnote{Corresponding author.}}\\
{\small  School of Mathematics and Statistics, Jiangsu Normal University,}\\ {\small  Xuzhou, Jiangsu 221116,
People's Republic
of China.}\\
{\small E-mails:luyong@jsnu.edu.cn, wuqimath@163.com}}

\date{}
\maketitle
\begin{center}
\begin{minipage}{120mm}
\vskip 0.3cm
\begin{center}
{\small {\bf Abstract}}
\end{center}
{\small A complex unit gain graph  is a triple
$\varphi=(G, \mathbb{T}, \varphi)$ (or $G^{\varphi}$ for short) consisting of a simple
graph $G$, as the underlying graph of $G^{\varphi}$, the set of unit
complex numbers $\mathbb{T}={z\in \mathbb{C}: |z| = 1}$ and a gain function
$\varphi: \overrightarrow{E}\rightarrow \mathbb{T}$ such that $\varphi(e_{i,j})=\varphi(e_{j,i}) ^{-1}$. Let $A(G^{\varphi})$ be adjacency matrix of $G^{\varphi}$. In this paper, we prove that
$$m(G)-c(G)\leq p(G^{\varphi})\leq m(G)+c(G),$$
$$m(G)-c(G)\leq n(G^{\varphi})\leq m(G)+c(G),$$
where $p(G^{\varphi})$, $n(G^{\varphi})$, $m(G)$ and $c(G)$ are the number of positive eigenvalues of $A(G^{\varphi})$, the number of negative eigenvalues of $A(G^{\varphi})$, the matching number and the cyclomatic number of $G$, respectively. Furthermore, we characterize
the graphs which attain the upper bounds and the lower
bounds, respectively.

 }
\end{minipage}
\end{center}

\section{Introduction }
 All graphs in this paper are simple graphs, without multiedges and loops. Let $G=(V(G), E(G))$
 be a simple graph and $V(G)=\{v_{1},v_{2},\ldots,v_{n}\}$, $E(G)$ be the vertex set and edge set of $G$, respectively. Denote by $e_{ij}$  the oriented edge from  $v_{i}$ to $v_{j}$. Let $\overrightarrow{E}$ be the oriented edge set obtained from $E(G)$, and $\{e_{ij},e_{ji}\}\in \overrightarrow{E}$. Denote by $N_{G}(x)$  the \emph{neighbor set} of a vertex $x\in V(G)$, and $d_{G}(x)=|N_{G}(x)|$  the \emph{degree} of $x$. In particular, $x$ is called a \emph{pendant} vertex of $G$ if $d_{G}(x)=1$. The neighbour vertex of a pendant vertex is called a \emph{quasi-pendant} vertex in graph $G$. We use $P_{n}$, $C_{n}$ to denote a \emph{path}, a \emph{cycle} of order $n$, respectively. An induced subgraph $C_{p}$ of a graph $G$ is called a \emph{pendant cycle} if $C_{p}$ is a cycle and has a unique vertex of degree 3 in $G$. Denote by $F(G)$ the edge set that every edge in $F(G)$  has an endpoint on some cycle of $G$ and other endpoint outside the cycle of $G$.

  Let $c(G)=|E(G)|-|V(G)|+\omega(G)$ be the \emph{cyclomatic number} of a graph $G$, where $\omega(G)$ is the number of connected components of $G$. Two distinct edges in a graph $G$ are \emph{independent} if they do not have common end-vertex in $G$. A set of pairwise independent edges of $G$ is called a \emph{matching} of $G$. A matching with the maximum cardinality is a \emph{maximum matching} of $G$, denote by $M(G)$. The \emph{matching number} of $G$, denoted by $m(G)$, is the cardinality of a maximum matching of $G$. A maximum matching $M(G)$ of $G$ \emph{saturates} $v$ of $G$ if $v$ is an endpoint of an edge in $M(G)$. Two graphs are called \emph{vertex-disjoint} if they have no common vertices.
  Let $G$ be a graph with pairwise vertex-disjoint cycles. $T_{G}$ is an acyclic graph obtained from $G$ by contracting each cycle of $G$ into a vertex, called a \emph{cyclic vertex}. Denote by $O(G)$  the set of vertices in cycles of $G$.

 Denote by $\mathfrak{G}$ a class of graphs that every graph $G\in\mathfrak{G}$ hold the following properties: (1) $G$ contains at least one cycle but  is not the disjoint union the disjoint cycles and/or trees, (2) any two cycles of $G$ share no common vertices if $G$ contains more than one cycle.

The \emph{adjacency matrix}  of $G$, denote by $A(G)$ is the symmetric $n\times n$ matrix with entries  $a_{ij}=1$ if and only if $v_{i}v_{j}\in E(G)$ and 0 elsewhere. Denote by $r(G)$ the \emph{rank} of $G$, which is the rank of $A(G)$. The multiplicity of 0 as an eigenvalue of $A(G)$, denoted by $\eta(G)$, is called the \emph{nullity} of $G$. The \emph{positive inertia index}(resp.the \emph{negative inertia index}), denoted by $p(G)$(resp.$n(G)$) is the number of positive eigenvalues(resp.negative eigenvalues) of $A(G)$. For a  simple graph of order $n$, $r(G)=p(G)+n(G)=n-\eta(G)$.

Denote by $\Phi=(G,\mathbb{T},\varphi)$ a complex unit gain graph, where $G$ is the \emph{underlying graph} of $\Phi$, $\mathbb{T}=\{z\in C: |z|=1\}$ is the \emph{circle group}, and   $\varphi:\overrightarrow{E}\rightarrow \mathbb{T}$, where $\varphi(e_{ij})=\varphi(e_{ji})^{-1}=\overline{\varphi(e_{ji})}$. For convenience, $\varphi(e_{ij})$ is also written as $\varphi_{v_{i}v_{j}}$ for $v_{i}v_{j}\in E(\Phi)$.  Denote by $G^{\varphi}$  a complex unit gain  graph $\Phi=(G,\mathbb{T},\varphi)$ for convenience. The \emph{adjacency matrix} of  $G^{\varphi}$ is the  Hermitian matrix $A(G^{\varphi})=(a_{ij})_{n\times n}$, where $a_{ij}=\varphi(e_{ij})$ if $v_{i}v_{j}\in E(G)$, and $a_{ij}=0$ otherwise. The rank $r(G^{\varphi})$ of $A(G^{\varphi})$ is called the \emph{rank}  of $G^{\varphi}$.The \emph{positive inertia index}(resp.the \emph{negative inertia index}), denoted by $p(G^{\varphi})$(resp.$n(G^{\varphi})$) is the number of positive eigenvalues(resp.negative eigenvalues) of $A(G^{\varphi})$. For a complex unit gain graph $G^{\varphi}$, the matching number, the cyclomatic number, pendant vertex and quasi-pendant vertex of $G^{\varphi}$ are defined to be the matching number, the cyclomatic number, pendant vertex and quasi-pendant vertex of its underlying graph, respectively.

 We write $G^{\varphi}-v$ ($v\in V(G^{\varphi})$) for the \emph{induced subgraph} obtained from $G^{\varphi}$ by deleting $v$ and all its incident edges. Let $V(G_{1}^{\varphi})\subset V(G^{\varphi})$, denote by $G^{\varphi}-G^{\varphi}_{1}$ the induced subgraph obtained from $G^{\varphi}$ by deleting all vertices of $G^{\varphi}_{1}$ and all incident edges. For $V(G^{\varphi}_{2})\subset V(G^{\varphi})$ and $v\notin G^{\varphi}_{2}$, denote by $G^{\varphi}_{2}+v$, the induced subgraph of $G^{\varphi}$ with vertex set $V(G^{\varphi}_{2})\cup \{v\}$. Let $G^{\varphi}$, $H^{\varphi}$ be two complex unit gain graphs, denote by $G^{\varphi}\cup H^{\varphi}$ the disjoin union of $G^{\varphi}$ and $H^{\varphi}$.

 Collatz et al.\cite{CS} had wanted to obtain all graphs of order $n$ with $r(G)<n$. Until today, this problem is also unsolved.
In mathematics, the rank (or nullity, inertia index) of a graph has attracted a lot of researchers' attention, they focus on the  bounds for  the rank (or nullity, inertia index) of a simple graph $G$ \cite{FYZW, FHLL, LGUO, RCZ, LWANG, WW, WLMA,ZWS}, a signed graph \cite{FDD,FWW, HHL,  LWZ1,LUWU1, YFQ}, an oriented graph \cite{LXYG, LWZ, LWZ2, MWTLAA, QHYG, WMT} and a mixed graph \cite{CHL,WLM} and so on.

For a complex unit gain graph $G^{\varphi}$, Reff \cite{NR}  defined the adjacency, incidence and Laplacian matrices of a complex unit gain graph. Some eigenvalue bounds for the adjacency and Laplacian matrices were present. Yu et al.\cite{YQT} give the inertia of some complex unit gain graph. Lu et al.\cite{LWX} characterized all the complex unit gain bicyclic graphs $G^{\varphi}$ with $r(G^{\varphi})=2,3,4$. Lu et al.\cite{LWZ3} obtained relation between the rank of a complex unit gain graph in terms of the rank of its underlying graph. Wang et al.\cite{WGF} obtained the  determinant of the Laplacian matrix of a complex unit gain graph. Xu et al.\cite{XZWT} characterized all the complex unit gain graphs of rank 2.  He et al.\cite{HHY} obtained the  bounds for the rank of a complex unit gain graph in terms of the independence number. Zaman and  He \cite{ZaHe} obtained the relation between the inertia indices of a complex unit gain graph and those of its underlying graph.  Lu and Wu \cite{LUWU} obtained the bounds for the rank of a complex unit gain graph in terms of its maximum degree.

In 2020, Li, Wang\cite{LiW} and He et al.\cite{HHD} obtained the rank of a complex unit graph in terms of the matching number, respectively. They proved that $$2m(G)-2c(G)\leq r(G^{\varphi})\leq 2m(G)+c(G).$$  All corresponding extremal graphs are characterized by them. Motivated by their results, in this paper, we will prove that $$m(G)-c(G)\leq p(G^{\varphi})\leq m(G)+c(G),$$
$$m(G)-c(G)\leq n(G^{\varphi})\leq m(G)+c(G).$$  All corresponding extremal graphs are characterized.

 In Section 2, we give and prove some lemmas about complex unit gain graphs. In Section 3, we characterize the relations between the inertia indices of a complex unit gain graph and its matching number.

\section{Preliminaries}
\quad For a simple graph $G$, there has some lemmas.
\noindent\begin{lemma}\label{le:2.1}\cite{CHL}
Let $G$ be a simple graph. Then $m(G)-1\leq m(G-v)\leq m(G)$ for any vertex $v\in V(G)$.
\end{lemma}
\noindent\begin{lemma}\label{le:2.2}\cite{HHD}
Let $G$ be a graph obtained by joining a vertex of an even cycle $C$ by an edge to a vertex of a connected graph $H$. Then $m(G)=m(C)+m(H)$.
\end{lemma}
\noindent\begin{lemma}\label{le:2.3}\cite{HHD}
Let $x$ be a pendant vertex of a graph G and $y$ be the neighbour of $x$. Then $m(G)=m(G-y)+1=m(G-x-y)+1$.
\end{lemma}
\noindent\begin{lemma}\label{le:2.4}\cite{FYZW}
Let $G$ be a graph with at least one cycle. Suppose that all cycles of $G$
are pairwise vertex-disjoint and each cycle is odd, then $m(T_{G})=m(G-O(G))$ if and only if there exists a maximum matching $M(G)$ of $G$ such that $M(G)\cap F(G)=\varnothing$.
\end{lemma}
\noindent\begin{lemma}\label{le:2.13}\cite{FYZW}
Let $G\in \mathfrak{G}$. If $m(T_{G})=m(G-O(G))$, then
$G$ contains at least one pendant vertex, and any quasi-pendant vertex of
$G$ lies outside of cycles.
\end{lemma}
\noindent\begin{lemma}\label{le:2.6}\cite{WMT}
Let $G$ be a graph with $x\in V (G)$. Then
\begin{enumerate}[(a)]

\item $c(G)=c(G-x)$ if $x$ lies outside any cycle of $G$;

\item $c(G-x)\leq c(G)-1$ if $x$ lies on a cycle of $G$;

\item $c(G-x)\leq c(G)-2$ if $x$ is a common vertex of distinct cycles of $G$.
\end{enumerate}
\end{lemma}
For a complex unit gain graph, we have the following definition and lemmas.

\noindent\begin{definition}\label{de:2.7}\cite{LWX}
Let $C_{n}^{\varphi} (n\geq3)$ be a complex unit gain cycle, denote by $$\varphi(C_{n})=\varphi_{v_{1}v_{2}}\varphi_{v_{2}v_{3}}\cdots \varphi_{v_{n-1}v_{n}}\varphi_{v_{n}v_{1}}.$$ Then $C_{n}^{\varphi}$ is said to be:

\begin{displaymath}
\left\{\
        \begin{array}{ll}
          \rm Type~A,&  \emph{if}~\varphi(C_{n})=(-1)^{n/2}~\emph{and}~n~\emph{is~even},\\
          \rm Type~B,& \emph{if}~\varphi(C_{n})\neq(-1)^{n/2}~\emph{and}~n~\emph{is~even},\\
          \rm Type~C,& \emph{if}~Re\left((-1)^{{(n-1)}/{2}}\varphi(C_{n})\right)>0~\emph{and}~n~\emph{is~odd},\\
          \rm Type~D,& \emph{if}~Re\left((-1)^{{(n-1)}/{2}}\varphi(C_{n})\right)<0~\emph{and}~n~\emph{is~odd},\\
          \rm Type~E,& \emph{if}~Re\left((-1)^{{(n-1)}/{2}}\varphi(C_{n})\right)=0~\emph{and}~n~\emph{is~odd}.
        \end{array}
      \right.
\end{displaymath}
\end{definition}
where $Re(\cdot)$ is the real part of a complex number.
\noindent\begin{lemma}\label{le:2.8}\cite{YQT}
Let $C_{n}^{\varphi}$ be a complex unit gain cycle of order $n$. Then

\begin{displaymath}
(p(C_{n}^{\varphi}),n(C_{n}^{\varphi}))=\left\{\
        \begin{array}{ll}
          (\frac{n-2}{2},\frac{n-2}{2}),& \emph{if}~C_{n}^{\varphi}\rm~is~of~Type~A,\\
          (\frac{n}{2},\frac{n}{2}),& \emph{if}~C_{n}^{\varphi}\rm~is~of~Type~B,\\
          (\frac{n+1}{2},\frac{n-1}{2}),& \emph{if}~C_{n}^{\varphi}\rm~is~of~Type~C,\\
          (\frac{n-1}{2},\frac{n+1}{2}),& \emph{if}~C_{n}^{\varphi}\rm~is~of~Type~D,\\
          (\frac{n-1}{2},\frac{n-1}{2}),& \emph{if}~C_{n}^{\varphi}\rm~is~of~Type~E.
        \end{array}
      \right.
\end{displaymath}
\end{lemma}

\noindent\begin{lemma}\label{le:2.9}\cite{YQT}
 Let $G^{\varphi}$ be a complex unit gain graph.
\begin{enumerate}[(a)]
\item Let $H^{\varphi}$ be an induced subgraph of $G^{\varphi}$. Then $p(H^{\varphi})\leq p(G^{\varphi})$ and $n(H^{\varphi})\leq n(G^{\varphi})$.
\item Let $G_{1}^{\varphi},G_{2}^{\varphi},\cdots,G_{t}^{\varphi}$ be the connected components of $G^{\varphi}$. Then $p(G^{\varphi})=\sum_{i=1}^{t} p(G_{i}^{\varphi})$ and $n(G^{\varphi})=\sum_{i=1}^{t} n(G_{i}^{\varphi})$.
\item $p(G^{\varphi})=0$ ($n(G^{\varphi})=0$) if and only if $G^{\varphi}$ is a graph without edges.
\end{enumerate}
\end{lemma}

\noindent\begin{lemma}\label{le:2.10}\cite{YQT}
Let $T^{\varphi}$ be an acyclic complex unit gain graph. Then $p(T^{\varphi})=n(T^{\varphi})=m(T)$.
\end{lemma}
\noindent\begin{lemma}\label{le:2.11}\cite{YQT}
Let $G^{\varphi}$ be a complex unit gain graph that has a pendant vertex $u$ and $v$ is the unique neighbour of $u$. Then $p(G^{\varphi})=p(G^{\varphi}-v)+1=p(G^{\varphi}-u-v)+1,n(G^{\varphi})=n(G^{\varphi}-v)+1=n(G^{\varphi}-u-v)+1$.
\end{lemma}
\noindent\begin{lemma}\label{le:2.12}\cite{ZaHe}
Let $G^{\varphi}$ be a complex unit gain graph and $u$ is a vertex of $G^{\varphi}$. Then $p(G^{\varphi})-1\leq p(G^{\varphi}-u)\leq p(G^{\varphi}), n(G^{\varphi})-1\leq n(G^{\varphi}-u)\leq n(G^{\varphi})$.
\end{lemma}

\section{Bounds for the inertia indices of $G^{\varphi}$}
\quad In this section, we will  obtain some bounds for the inertia indices of a complex unit gain graph in terms of its matching number.

\noindent\begin{theorem}\label{th:3.1}
Let $G^{\varphi}$ be a  connected complex unit gain graph. Then
\begin{enumerate}[(a)]
\item $m(G)-c(G)\leq p(G^{\varphi})\leq m(G)+c(G)$.
\item $m(G)-c(G)\leq n(G^{\varphi})\leq m(G)+c(G)$.
\end{enumerate}
\end{theorem}
\noindent\textbf{Proof.}
we get the inequalities above by induction on $c(G)$. If $c(G)=0$, then $G^{\varphi}$ is a tree, $p(G^{\varphi})=n(G^{\varphi})=m(G)$ holds by Lemma \ref{le:2.10}.
Now we assume that assertion holds for every connected complex unit gain graph that the cyclomatic number is less than $c(G)$. Let $v$ be a vertex on some cycle in $G^{\varphi}$. Let $H_{1},\cdots,H_{k}$ be the connected components of $G-v$. By Lemma \ref{le:2.6}(b), we have
$$\sum_{i=1}^{k} c(H_{i})=c(G-v)\leq c(G)-1.$$
Combining with  Lemmas \ref{le:2.1}, \ref{le:2.9}(b) and \ref{le:2.12}, we have
$$\sum_{i=1}^{k} p(H_{i}^{\varphi})=p(G^{\varphi}-v)\leq p(G^{\varphi})\leq p(G^{\varphi}-v)+1$$
and
$$m(G)-1\leq m(G-v)=\sum_{i=1}^{k} m(H_{i})\leq m(G).$$

Because $c(G-v)\leq c(G)-1$, by the induction hypothesis, for each $i\in\{1,2,\ldots,k\}$,
$$m(H_{i})-c(H_{i})\leq p(H_{i}^{\varphi})\leq m(H_{i})+c(H_{i})$$

Hence,
\begin{align*}
&p(G^{\varphi})\geq p(G^{\varphi}-v)\\
&~~~~~~~~=\sum_{i=1}^{k} p(H_{i}^{\varphi})\\
&~~~~~~~~\geq \sum_{i=1}^{k} [m(H_{i})-c(H_{i})]\\
&~~~~~~~~\geq m(G)-1-(c(G)-1)\\
&~~~~~~~~=m(G)-c(G).
\end{align*}
and
\begin{align*}
&p(G^{\varphi})\leq p(G^{\varphi}-v)+1\\
&~~~~~~~~=\sum_{i=1}^{k} p(H_{i}^{\varphi})+1\\
&~~~~~~~~~\leq \sum_{i=1}^{k} [m(H_{i})+c(H_{i})]+1\\
&~~~~~~~~~\leq m(G)+c(G)-1+1\\
&~~~~~~~~~=m(G)+c(G).
\end{align*}
That is $$m(G)-c(G)\leq p(G^{\varphi})\leq m(G)+c(G).$$
Similarly, we can obtain that
$$m(G)-c(G)\leq n(G^{\varphi})\leq m(G)+c(G).$$

This completes the proof.  \quad $\square$

For convenience, we call $G^{\varphi}$ to be $p$-lower optimal(resp.,$p$-upper optimal) if $p(G^{\varphi})$ obtains the lower bound (upper bound) in Theorem \ref{th:3.1}. Similarly, $G^{\varphi}$ is $n$-lower optimal(resp.,$n$-upper optimal) if $n(G^{\varphi})$ obtains the lower bound (upper bound) in Theorem \ref{th:3.1}.

\noindent\begin{lemma}\label{le:3.3}
Let $G^{\varphi}$ be a complex unit gain graph and $H_{1}^{\varphi},H_{2}^{\varphi},\cdots,H_{k}^{\varphi}$ are its connected components. Then
\begin{enumerate}[(a)]
  \item $G^{\varphi}$ is $p$-lower optimal if and only if for each $i\in\{1,2,\ldots,k\}$, $H_{i}^{\varphi}$ is $p$-lower optimal.
  \item $G^{\varphi}$ is $p$-upper optimal if and only if for each $i\in\{1,2,\ldots,k\}$, $H_{i}^{\varphi}$ is $p$-upper optimal.
\end{enumerate}

\end{lemma}
\noindent\textbf{Proof.} At first, we will prove the (a) of this lemma.

\textbf{Necessity:} Assume that the conclusion is not true, without loss of generality, we suppose that $H_{1}^{\varphi}$ is not $p$-lower optimal. By Theorem \ref{th:3.1},
$$p(H_{1}^{\varphi})>m(H_{1})-c(H_{1}),$$
and for each $i\in\{2,3,\ldots,k\}$, we have
$$p(H_{i}^{\varphi})\geq m(H_{i})-c(H_{i}).$$

By Lemma \ref{le:2.9}(b), we have $$p(G^{\varphi})=\sum_{i=1}^{k} p(H_{i}^{\varphi})>\sum_{i=1}^{k} [m(H_{i})-c(H_{i})]=m(G)-c(G),$$
a contradiction.

\textbf{Sufficiency:} For each $i\in\{1,2,\ldots,k\}$, $H_{i}^{\varphi}$ is $p$-lower optimal, so we have
$$p(H_{i}^{\varphi})=m(H_{i})-c(H_{i}).$$

By Lemma \ref{le:2.9}(b), we have $$p(G^{\varphi})=\sum_{i=1}^{k} p(H_{i}^{\varphi})=\sum_{i=1}^{k} [m(H_{i})-c(H_{i})]=m(G)-c(G).$$

Using the same method, we can prove the (b) of this lemma.

This complete the proof.   \quad $\square$

\noindent\begin{lemma}\label{le:3.4}
Let $G^{\varphi}$ be a complex unit gain graph and a vertex $v$ lies on some cycle of $G^{\varphi}$.
Then the following results are established:
\begin{enumerate}[(a)]
\item If $G^{\varphi}$ is $p$-lower optimal, then $p(G^{\varphi})=p(G^{\varphi}-v)$, $p(G^{\varphi}-v)=m(G-v)-c(G-v)$, $c(G)=c(G-v)+1$, $m(G)=m(G-v)+1$;
\item If $G^{\varphi}$ is $p$-upper optimal, then $p(G^{\varphi})=p(G^{\varphi}-v)+1$, $p(G^{\varphi}-v)=m(G-v)+c(G-v)$, $c(G)=c(G-v)+1$, $m(G)= m(G-v)$;
\item If $G^{\varphi}$ is $p$-lower optimal (or upper optimal), then $v$ lies on just one cycle of $G^{\varphi}$ and $v$ is not a quasi-pendant vertex.
\end{enumerate}
\end{lemma}
\noindent\textbf{Proof.}
By the proof of Theorem \ref{th:3.1}, we can obtain the (a) and (b). For (c), when $G^{\varphi}$ is $p$-lower optimal. If $v$ lie on at least two cycles of $G^{\varphi}$, by Lemma \ref{le:2.6}(c),
$$c(G-v)\leq c(G)-2,$$ which contradicts (a). If $v$ is a quasi-pendant vertex of $G^{\varphi}$, then by Lemma \ref{le:2.11},
$$p(G^{\varphi}-v)=p(G^{\varphi})-1,$$ which contradicts (a).

When $G^{\varphi}$ is $p$-upper optimal. If $v$ lie on at least two cycles of $G^{\varphi}$, by Lemma \ref{le:2.6}(c),
$$c(G-v)\leq c(G)-2,$$ which contradicts (b). If $v$ is a quasi-pendant vertex of $G^{\varphi}$, then by Lemma \ref{le:2.3},
$$m(G)=m(G-v)+1,$$ which contradicts (b).
 So the assertion (c) is  hold.   \quad $\square$

\noindent\begin{lemma}\label{le:3.5}
Let  $G^{\varphi}$ be a complex unit gain graph which contains a pendent vertex $x$ with its unique neighbour $y$. Let $H^{\varphi}=G^{\varphi}-x-y$. If $G^{\varphi}$ is $p$-lower optimal (upper optimal), then $H^{\varphi}$ is also $p$-lower optimal (upper optimal).
\end{lemma}
\noindent\textbf{Proof.}
If $G^{\varphi}$ is $p$-lower optimal, we have
$$p(G^{\varphi})=m(G)-c(G).$$

By Lemmas \ref{le:2.3}, \ref{le:2.6}(a) and \ref{le:2.11}, we have $$p(H^{\varphi})=p(G^{\varphi})-1=m(G)-c(G)-1=m(H)+1-c(H)-1=m(H)-c(H),$$
so $H^{\varphi}$ is also $p$-lower optimal.

Using the same method, we can get that if $G^{\varphi}$ is $p$-upper optimal, then $H^{\varphi}$ is $p$-upper optimal.
  \quad $\square$

\noindent\begin{lemma}\label{le:3.6}
Let $G^{\varphi}$ be a complex unit gain unicyclic graph and $C_{q}^{\varphi}$ is unique  cycle of $G^{\varphi}$.
\begin{enumerate}[(a)]
  \item If $G^{\varphi}$ is $p$-lower optimal, then $C_{q}^{\varphi}$ is of Type A.
  \item If $G^{\varphi}$ is $p$-upper optimal, then $C_{q}^{\varphi}$ is of Type C.
\end{enumerate}
\end{lemma}
\noindent\textbf{Proof.} At first, we will prove the (a) of this lemma. We shall apply induction on the order of $T_{G}$.
 If $|V(T_{G})|=1$, then $G^{\varphi}$ is $C_{q}^{\varphi}$.  By Lemma \ref{le:2.8}, we can get that $C_{q}^{\varphi}$ is of Type A.
 Now suppose that $|V(T_{G})|\geq2$, if $G^{\varphi}$ has no pendent vertices, $G^{\varphi}$ is the union of $C_{q}^{\varphi}$ and some isolated vertices, by Lemmas \ref{le:2.8} and \ref{le:2.9}(b), $C_{q}^{\varphi}$ is of Type A. If $G^{\varphi}$ has a pendent vertex, say $u$ and $v$ is its neighbour vertex in $G^{\varphi}$. By Lemma \ref{le:3.4}(c), $v$ is not on $C_{q}^{\varphi}$, so $|V(T_{G})|\neq 2$, $|V(T_{G})|\geq 3$. 
Then by Lemma \ref{le:3.5}, $G^{\varphi}-u-v$ is $p$-lower optimal and $C_{q}^{\varphi}$ is its unique cycle, by induction hypothesis, $C_{q}^{\varphi}$ is of Type A.

Using the same method, we can obtain the (b) of this lemma.
\quad $\square$

\noindent\begin{lemma}\label{le:3.7}
Let $H^{\varphi}$ be a complex unit gain graph with any two cycles (if any) share no common vertices. Denote by $G^{\varphi}$ the graph obtained from adding an edge between a vertex $u$ of a  cycle $C_{q}^{\varphi}$ and a vertex $v$ of $H^{\varphi}$. If $G^{\varphi}$ is $p$-lower optimal, then
\begin{enumerate}[(a)]
\item Every  cycle of $G^{\varphi}$ is of Type $A$;
\item The edge $uv$ does not belong to any maximum matching of $G$;
\item Each maximum matching of $H$ saturates $v$;
\item $m(H+u)=m(H)$; 
\item $H^{\varphi}$ is $p$-lower optimal;
\item Let $K^{\varphi}$ be the induced  subgraph of $G^{\varphi}$ with vertex set
$V(H)\cup{x}$. Then $K^{\varphi}$ is also $p$-lower optimal;
\end{enumerate}
\end{lemma}
\noindent\textbf{Proof.}
We use induction on $c(G)$, since $G^{\varphi}$ has a cycle $C_{q}^{\varphi}$, $c(G)\geq 1$.
If $H^{\varphi}$ contains no cycle, by Lemma \ref{le:3.6}(a), we have that $C_{q}^{\varphi}$ is of Type A.

If $H^{\varphi}$ contains at least one cycle. Let $x$ be a vertex lying on some cycle of $H^{\varphi}$ and $G_{0}^{\varphi}=G^{\varphi}-x$. By Lemma \ref{le:3.4}(a), we have $G_{0}^{\varphi}$ is $p$-lower optimal.

By induction hypothesis, one has that each cycle in $G_{0}^{\varphi}$, including $C_{q}^{\varphi}$ is of Type A. By a
similar discussion as for $G^{\varphi}-u$, we can show that each cycle in $H^{\varphi}$ is of Type A.
This completes the proof of (a).

Suppose on the contrary that $uv$ belongs to a maximum matching
$M$ of $G$. By (a),  $C_{q}^{\varphi}$ is an even cycle, so there exists a vertex $w$ lying on $C_{q}^{\varphi}$
such that $w$ is not saturated by $M$. Thus we have
$$m(G)=m(G-w),$$ a contradiction to Lemma \ref{le:3.4}(a). This completes the proof of (b).

By Lemma \ref{le:2.2}, we have $$m(G)=m(H)+m(C_{q}).$$ Let $M_{1}$ be the maximum
matching of $C_{q}$. Suppose on the contrary that there exists a maximum matching $M_{2}$ of $H$ fails to saturate $v$, so $M_{1} \cup M_{2}$ is a maximum matching of $G$. Then we obtain a maximum matching $M^{'} \cup M_{2}$ of $G$ which contains $uv$, where $M^{'}$ is
obtained from $M_{1}$ by replacing the edge in $M_{1}$ which saturates $u$ with $uv$, a contradiction
to (b). This completes the proof of (c).

By (c), we can obtain the (d) of this lemma.

 By Lemma \ref{le:3.4}(a), $G^{\varphi}-u$ is $p$-lower optimal. Then (e) immediately follows
from Lemma \ref{le:3.3}(a).

Suppose that $C_{q}=uu_{2}u_{3} \cdots u_{2s}u$. Since $G^{\varphi}$ is $p$-lower optimal, by Lemma \ref{le:3.4}(a), one has that $G_{1}^{\varphi}=G^{\varphi}-u_{2}$ is also $p$-lower optimal. Obviously, $u_{3}$ and $u_{4}$ are
pendant vertex and quasi-pendant vertex of $G_{1}^{\varphi}$, respectively. By Lemma \ref{le:3.5}, one has
that $G_{2}^{\varphi}=G_{1}^{\varphi}-u_{3}-u_{4}$ is also $p$-lower optimal. Repeating such process (deleting
a pendant vertex and a quasi-pendant vertex), after $s-1$ steps, the result graph is
$$G^{\varphi}-u_{2}-u_{3}-\cdots-u_{2s}=H^{\varphi}+u=K^{\varphi}.$$ By Lemma \ref{le:3.5}, $K^{\varphi}$ is also $p$-lower optimal.
   \quad $\square$

\noindent\begin{lemma}\label{le:3.8}
Let $G^{\varphi}$ be a connected complex unit gain graph. If $G^{\varphi}$ is $p$-lower optimal,
then there exists a maximum matching $M$ of $G$ such that $M \cap F(G)=\varnothing$. Moreover, $m(G)=m(G-O(G))+\sum_{C\subseteq G} m(C)$, where $C$ goes through all cycles of $G$.
\end{lemma}
\noindent\textbf{Proof.}
We shall apply induction on the order of $G^{\varphi}$. If $G^{\varphi}$ is an isolated vertex, the result holds trivially. Suppose $G^{\varphi}$ contains at least two vertices.
If $G^{\varphi}$ has a pendant vertex $u$, and $v$ is the unique neighbor of $u$. Let $H^{\varphi}=G^{\varphi}-u-v$
and $H_{1}^{\varphi},H_{2}^{\varphi},\cdots,H_{t}^{\varphi}$ be all connected components of $H^{\varphi}$. By Lemmas \ref{le:3.4}(c) and \ref{le:3.5}, $v$ is not on the cycle of $G^{\varphi}$ and $H^{\varphi}$ is $p$-lower optimal. Then by Lemma \ref{le:3.3}(a), we have $H_{i}^{\varphi}$ is also $p$-lower optimal for each $i\in\{1,2,\ldots,t\}$.

By induction hypothesis, for $H_{i}^{\varphi}$ there exists a maximum matching $M_{i}$ of $H_{i}$ such that $M_{i}\cap F(H_{i})=\varnothing$ for each $i\in \{1,2,\ldots,t\}$. Let $$M=(\cup_{i=1}^{t} M_{i} )\cup\{uv\}.$$ We can obtain  that $M$ is a
maximum matching of $G$ which satisfies $$M\cap F(G)=\varnothing.$$

If $G^{\varphi}$ has no pendant vertices, then $G^{\varphi}$ contains a pendant complex unit gain cycle, says $C^{\varphi}$. Let $$G_{1}^{\varphi}=G^{\varphi}-C^{\varphi}. $$ By Lemmas \ref{le:3.7}(a) and (e), it follows that $G_{1}^{\varphi}$ is $p$-lower optimal and the order of cycle $C^{\varphi}$ is even. Applying the induction on $G_{1}^{\varphi}$, there exists a maximum matching $M_{0}$ of $G_{1}$ such that $$M_{0} \cap F(G_{1})=\varnothing.$$

Let $M_{1}$ be a maximum matching of $C$. By Lemma \ref{le:2.2},  $M=M_{0}\cup M_{1}$
is a maximum matching of $G$ satisfying $$M\cap F(G)=\varnothing.$$  Moreover, we can get that
$$m(G)=m(G-O(G))+\sum_{C\subseteq G}m(C),$$ where $C$ goes through all cycles of $G$.

This completes the proof.   \quad $\square$
\noindent\begin{theorem}\label{th:3.9}
Let $G^{\varphi}$ be a connected complex unit gain graph. Then $G^{\varphi}$ is $p$-lower optimal if and only if the following
three conditions all hold:
\begin{enumerate}[(a)]
\item Any two cycles of $G^{\varphi}$  share no common vertices;
\item Each cycle of $G^{\varphi}$ is of Type A;
\item $m(T_{G})= m(G-O(G))$.
\end{enumerate}
\end{theorem}
\noindent\textbf{Proof.}
\textbf{Sufficiency:} We shall apply induction on the order of $G^{\varphi}$. If $G^{\varphi}$ is a complex unit gain
tree or single cycle of Type A,  the result follows from Lemmas \ref{le:2.8} and
\ref{le:2.10}. If not, since $m(T_{G})=m(G-O(G))$, we have that $G^{\varphi}$ contains at least one pendant vertex, say $x$. Let $y$ be the unique neighbour vertex of $x$ in $G^{\varphi}$. By Lemma \ref{le:2.13}, $y$ lies outside any cycle of $G^{\varphi}$.

Let $H^{\varphi}=G^{\varphi}-x-y$ and $H_{1}^{\varphi},H_{2}^{\varphi},\cdots,H_{t}^{\varphi}$ be all connected components
of $H^{\varphi}$.

Since $x,y$ are both outside the cycle of $G^{\varphi}$, then $$T_{H}=T_{G}-x-y$$ and
$$H-O(H)=G-O(G)-x-y.$$ By condition(c) and Lemma \ref{le:2.3}, we have
\begin{align*}
&m(T_{G})=\sum_{i=1}^{t} m(T_{H_{i}})+1\\
 &~~~~~~~~~=m(G-O(G))\\
 &~~~~~~~~~=\sum_{i=1}^{t} m(H_{i}-O(H_{i}))+1.
 \end{align*}

Therefore $$\sum_{i=1}^{t} m(T_{H_{i}})=\sum_{i=1}^{t} m(H_{i}-O(H_{i})).$$ We observe that for each $i\in\{1,2,\ldots,t\}$, $$m(T_{H_{i}})\geq m(H_{i}-O(H_{i})),$$ so we obtain $$m(T_{H_{i}})=m(H_{i}-O(H_{i}))$$ for each $i\in\{1,2,\ldots,t\}$. Therefore, $H_{i}^{\varphi}$ satisfies (a)-(c) for each $i\in \{1,2,\ldots,t\}$.

By induction hypothesis, for each $i\in \{1,2,\ldots,t\}$,
$$p(H_{i}^{\varphi})=m(H_{i})-c(H_{i}).$$
Then by Lemmas \ref{le:2.3}, \ref{le:2.9} and \ref{le:2.11},
\begin{align*}
&p(G^{\varphi})=p(H^{\varphi})+1\\
 &~~~~~~~~=\sum_{i=1}^{t} p(H_{i}^{\varphi})+1\\
 &~~~~~~~~=\sum_{i=1}^{t} [m(H_{i})-c(H_{i})]+1\\
 &~~~~~~~~=m(H)-c(H)+1\\
 &~~~~~~~~=m(G)-c(G).
 \end{align*}

\textbf{Necessity:} Let $G^{\varphi}$ be a complex unit gain graph such that $$p(G^{\varphi})=m(G)-c(G).$$ If $G^{\varphi}$ is a complex unit gain tree, $G^{\varphi}$ clearly satisfies (a)-(c) of this theorem. Assume that $G^{\varphi}$ has at least one complex unit gain cycle. The
assertion (a) follows from Lemma \ref{le:3.4}(c).

For (b), if $c(G)=1$, the result holds by Lemma \ref{le:3.6}. Now assume $c(G)=k$, where $k\geq 2$. If there exists a cycle, say $C_{1}^{\varphi}$, which is not of Type A, then by deleting an arbitrary vertex of each complex unit gain cycle of $G^{\varphi}$ except for
$C_{1}^{\varphi}$, we get a complex unit gain graph $H^{\varphi}$ with $c(H)=1$ and by Lemma \ref{le:3.6}, $$p(H^{\varphi})\geq m(H).$$ By Lemmas \ref{le:2.1} and \ref{le:2.9}(a),
\begin{align*}
&p(G^{\varphi})\geq p(H^{\varphi})\\
 &~~~~~~~~\geq m(H)\\
 &~~~~~~~~\geq m(G)-(k-1)\\
 &~~~~~~~~=m(G)-c(G)+1,
\end{align*}
a contradiction.

We prove the assertion (c) by the induction on the order of $G^{\varphi}$. If $G^{\varphi}$ is a complex unit gain cycle or tree, the result follows. If not, we can consider the following two cases.

\textbf{Case 1.} $G^{\varphi}$ contains a pendant vertex, say $x$.
Let $y$ be the unique neighbour vertex of $x$ in $G^{\varphi}$ and $$H^{\varphi}=G^{\varphi}-x-y.$$ Let $H_{1}^{\varphi},H_{2}^{\varphi},\cdots,H_{t}^{\varphi}$ be all connected components of $H^{\varphi}$. By Lemmas \ref{le:3.4}(c)
and \ref{le:3.5}, $y$ does not lie on any cycle of $G^{\varphi}$ and $H^{\varphi}$ is also $p$-lower optimal. Then
by Lemma \ref{le:3.3}(a), $H_{i}^{\varphi}$ is $p$-lower optimal for each $i\in \{1,2,\ldots,t\}$.

By induction hypothesis, we have $m(T_{H_{i}})=m(H_{i}-O(H_{i}))$ for each $i\in \{1,2,\ldots,t\}$.
Then by Lemma \ref{le:2.3}, we have that
\begin{align*}
&m(T_{G})=m(T_{H})+1\\
 &~~~~~~~~~=\sum_{i=1}^{t}m(T_{H_{i}})+1\\
 &~~~~~~~~~=\sum_{i=1}^{t}m(H_{i}-O(H_{i}))+1\\
 &~~~~~~~~~=m(H-O(H))+1\\
 &~~~~~~~~~=m(G-O(G)).
\end{align*}
Assertion (c) holds in this case.

\textbf{Case 2.} $G^{\varphi}$ has a pendant complex unit gain cycle, say $C_{q_{1}}^{\varphi}$. Let $x$ be the unique vertex with degree $3$ in $C_{q_{1}}^{\varphi}$. Let $C_{q_{1}}^{\varphi},C_{q_{2}}^{\varphi},\cdots,C_{q_{t}}^{\varphi}$ be all cycles of $G^{\varphi}$ and
$K^{\varphi}=G^{\varphi}-C_{q_{1}}^{\varphi},H^{\varphi}=K^{\varphi}+x$. By Lemma \ref{le:3.7}(f), one has that $H^{\varphi}$ is $p$-lower optimal. Because $|V(H^{\varphi})|< |V(G^{\varphi})|$, by induction hypothesis, we have
$$m(T_{H})=m(H-O(H)).$$
Note that $T_{G}\cong T_{H}$. So by Lemmas \ref{le:3.7}(d) and \ref{le:3.8}, we have
\begin{align*}
&m(T_{G})=m(T_{H})\\
 &~~~~~~~~~=m(H-O(H))\\
 &~~~~~~~~~=m(H)-\frac{\sum_{i=2}^{t}|V(C_{q_{i}})|}{2}\\
 &~~~~~~~~~=m(K)-\frac{\sum_{i=2}^{t}|V(C_{q_{i}})|}{2}\\
 &~~~~~~~~~=m(K)+\frac{|V(C_{q_{1}})|}{2}-\frac{\sum_{i=1}^{t}|V(C_{q_{i}})|}{2}\\
 &~~~~~~~~~=m(G)-\frac{\sum_{i=1}^{t}|V(C_{q_{i}})|}{2}\\
 &~~~~~~~~~=m(G-O(G)).
\end{align*}

This completes the proof.   \quad $\square$


\noindent\begin{theorem}\label{th:3.10}
Let $G^{\varphi}$ be a connected complex unit gain graph. Then $G^{\varphi}$ is $p$-upper optimal if and only if the following
three conditions all hold:
\begin{enumerate}[(a)]
 \item Any two cycles of $G^{\varphi}$ share no common vertices;
 \item Each cycle of $G^{\varphi}$ is of Type C;
 \item $m(T_{G})=m(G-O(G))$.
\end{enumerate}
\end{theorem}
\noindent\textbf{Proof.} \textbf{Sufficiency:} We will use induction on the order of $G^{\varphi}$. If $G$ is a complex unit gain
tree or cycle of Type C, the result follows from Lemmas \ref{le:2.8} and \ref{le:2.10}. If not, since $$m(T_{G})=m(G-O(G)),$$ by Lemma \ref{le:2.13}, we have that $G^{\varphi}$ contains at least one pendant vertex, say $x$. Let $y$ be the unique neighbour vertex of $x$ in $G^{\varphi}$. By Lemma \ref{le:2.13}, $y$ lies outside
any cycle of $G^{\varphi}$. Let $$H^{\varphi}=G^{\varphi}-x-y$$ and $H_{1}^{\varphi},H_{2}^{\varphi},\cdots,H_{t}^{\varphi}$ be all connected components of $H^{\varphi}$. Since $x,y$ are both outside the cycle of $G^{\varphi}$, then $$T_{H}=T_{G}-x-y$$ and
$$H-O(H)=G-O(G)-x-y.$$ By condition(c) and Lemma \ref{le:2.3}, we have
\begin{align*}
&m(T_{G})=\sum_{i=1}^{t} m(T_{H_{i}})+1\\
&~~~~~~~~~=m(G-O(G))\\
&~~~~~~~~~=\sum_{i=1}^{t} m(H_{i}-O(H_{i}))+1.
\end{align*}

Therefore $$\sum_{i=1}^{t} m(T_{H_{i}})=\sum_{i=1}^{t} m(H_{i}-O(H_{i})),$$ we observe that for each $i\in \{1,2,\ldots,t\}$, $$m(T_{H_{i}})\geq m(H_{i}-O(H_{i})),$$ so we obtain $$m(T_{H_{i}})=m(H_{i}-O(H_{i}))$$ for each $i\in \{1,2,\ldots,t\}$. Therefore, $H_{i}^{\varphi}$ satisfies (a)-(c) for each $i\in \{1,2,\ldots,t\}$.

By induction hypothesis, for each $i\in\{1,2,\ldots,t\}$,
$$p(H_{i}^{\varphi})=m(H_{i})+c(H_{i}).$$
Then by Lemmas \ref{le:2.3}, \ref{le:2.9} and \ref{le:2.11},
\begin{align*}
&p(G^{\varphi})=p(H^{\varphi})+1\\
 &~~~~~~~~=\sum_{i=1}^{t} p(H_{i}^{\varphi})+1\\
 &~~~~~~~~=\sum_{i=1}^{t} [m(H_{i})+c(H_{i})]+1\\
 &~~~~~~~~=m(H)+c(H)+1\\
 &~~~~~~~~=m(G)+c(G).
\end{align*}

\textbf{Necessity:} Let $G^{\varphi}$ be a complex unit gain graph such that $$p(G^{\varphi})=m(G)+c(G).$$ The proof for $(a)$ and
$(b)$ goes parallel as in Theorem \ref{th:3.9}, thus omitted.
We prove the assertion $(c)$ by the induction on the order of $G^{\varphi}$. If $G^{\varphi}$ is a complex unit gain tree or cycle of Type C, the result follows. If not, we can consider the following two  cases.

\textbf{Case 1.} $G^{\varphi}$ contains a pendant vertex, say $x$.
Let $y$ be the unique neighbour vertex of $x$ in $G^{\varphi}$ and $H^{\varphi}=G^{\varphi}-x-y$. Let $H_{1}^{\varphi},H_{2}^{\varphi},\cdots,H_{t}^{\varphi}$ be all connected components of $H^{\varphi}$. By Lemmas \ref{le:3.4}(c)
and \ref{le:3.5}, $y$ does not lie on any cycle of $G^{\varphi}$ and $H^{\varphi}$ is also $p$-upper optimal. Then
by Lemma \ref{le:3.3}(b), $H_{i}^{\varphi}$ is $p$-upper optimal for each $i\in\{1,2,\ldots,t\}$.
By induction hypothesis, we have
$$m(T_{H_{i}})=m(H_{i}-O(H_{i}))$$ for each $i\in \{1,2,\ldots,t\}$.
Then by Lemma \ref{le:2.3},
\begin{align*}
&m(T_{G})=m(T_{H})+1\\
&~~~~~~~~~=\sum_{i=1}^{t} m(T_{H_{i}})+1\\
&~~~~~~~~~=\sum_{i=1}^{t} m(H_{i}-O(H_{i}))+1\\
&~~~~~~~~~=m(H-O(H))+1\\
&~~~~~~~~~=m(G-O(G)).
\end{align*}

Assertion (c) holds in this case.

\textbf{Case 2.} $G^{\varphi}$ has a pendant complex unit gain cycle, say $C_{p_{1}}^{\varphi}$. Let $x$ be the unique vertex with degree $3$ in $C_{p_{1}}^{\varphi}$. Let $C_{p_{1}}^{\varphi},C_{p_{2}}^{\varphi},\cdots,C_{p_{t}}^{\varphi}$ be all cycles of $G^{\varphi}$ and $H^{\varphi}=G^{\varphi}-C_{p_{1}}^{\varphi},K^{\varphi}=H^{\varphi}+x$, where $H^{\varphi}$ is a connected component of $G^{\varphi}-x$. By Lemma \ref{le:3.4}(b), one has that
$$p(G^{\varphi}-x)=m(G-x)+c(G-x)$$
and
$$m(G-x)=m(G).$$ Then by Lemma \ref{le:3.3}(b), $H^{\varphi}$ is $p$-upper optimal. Since $C_{p_{1}}$ is odd cycle and
$m(G-x)=m(G)$, we have that $$m(G)=m(C_{p_{1}})+m(H).$$ Because $|V(H^{\varphi})|< |V(G^{\varphi})|$, by induction hypothesis, we have
$$m(T_{H})=m(H-O(H)).$$

So by Lemma \ref{le:2.4}, there exists a maximum matching $M(H)$ of $H$ such that $M(H)\cap F(H)=\emptyset$. Let $M(C_{p_{1}})$
be a maximum matching of $C_{p_{1}}$. Then $$M(G)=M(H)\cup M(C_{p_{1}})$$ is a maximum matching of
$G$, which satisfies $$M(G)\cap F(G)=\emptyset.$$ Again by Lemma \ref{le:2.4}, we get $$m(T_{G})=m(G-O(G)).$$

This completes the proof.   \quad $\square$

Using the same methods as in Theorems \ref{th:3.9} and \ref{th:3.10}, we can obtain the following two theorems.
\noindent\begin{theorem}\label{th:3.11}
Let $G^{\varphi}$ be a connected complex unit gain graph. Then $G^{\varphi}$ is $n$-lower optimal if and only if the following
three conditions all hold:
\begin{enumerate}[(a)]
\item Any two cycles of $G^{\varphi}$ share no common vertices;
\item Each cycle of $G^{\varphi}$ is of Type A;
\item $m(T_{G})= m(G-O(G))$.
\end{enumerate}
\end{theorem}

\noindent\begin{theorem}\label{th:3.12}
Let $G^{\varphi}$ be a connected complex unit gain graph. Then $G^{\varphi}$ is $n$-upper optimal if and only if the following
three conditions all hold:
\begin{enumerate}[(a)]
\item Any two cycles of $G^{\varphi}$ share no common vertices;
\item Each cycle of $G^{\varphi}$ is of Type D;
\item $m(T_{G})= m(G-O(G))$.
\end{enumerate}
\end{theorem}
\noindent\begin{remark}
 For a  complex unit gain graph $G^{\varphi}$, if $\varphi(\overrightarrow{E})\subset\{1\}$, then $G^{\varphi}$ is the underlying  graph $G$. If $\varphi(\overrightarrow{E})\subset\{1,-1\}$, then $G^{\varphi}$ is the signed graph $\Gamma$. If $\varphi(\overrightarrow{E})\subset\{1,i,-i\}$, then $G^{\varphi}$ is the mixed graph $D_{G}$. Combing with above, we know that the  results of complex unit gain graphs are also applies to simple graphs, signed graphs and mixed graphs.
\end{remark}

\end{document}